\documentclass[letterpaper, 10 pt, conference]{ieeeconf}  

\IEEEoverridecommandlockouts                              

\usepackage{graphics} 
\usepackage{epsfig} 
\usepackage{amsmath,amssymb,amsfonts} 
\usepackage{textcomp}
\usepackage{hyperref}
\usepackage{fancyhdr}
\usepackage{colortbl}
\usepackage{xcolor}
\usepackage{tikz}
\usepackage{graphicx,float}
\newtheorem{Theorem}{Theorem}
\newtheorem{Corollary}[Theorem]{Corollary}
\newtheorem{Proposition}[Theorem]{Proposition}

\newtheorem{Definition}[Theorem]{Definition}
\newtheorem{Remark}[Theorem]{Remark}
\usepackage{color} 

\title{\LARGE \bf
Stability and optimality of multi-scale transportation networks with distributed dynamic tolls
}

\author{Rosario Maggistro$^{1}$ and Giacomo Como$^{1,2}$
\thanks{*This reasearch was carried on within the framework of the MIUR-funded {\it Progetto di Eccellenza} of the {\it Dipartimento di Scienze Matematiche G.L.~Lagrange}, CUP: E11G18000350001, and was partly supported by the {\it Compagnia di San Paolo} and the Swedish Research Council.}
\thanks{$^{1}$Department of Mathematical Sciences, Politecnico di Torino,
	Corso Duca degli Abruzzi 24, 10129 Torino, Italy {\tt\small \{rosario.maggistro, giacomo.como\}@polito.it}.} 
\thanks{$^{2}$ Department of Automatic Control,
	Lund University, BOX 118, SE-22100 Lund, Sweden
	{\tt\small giacomo.como@control.lth.se}.} %
}

\begin{document}

\maketitle
\thispagestyle{empty}
\pagestyle{empty}
\begin{abstract}

We study transportation networks controlled by dynamical feedback tolls. We consider a multiscale transportation network model whereby the dynamics of the traffic flows are intertwined with those of the drivers' route choices. The latter are influenced by the congestion status on the whole network as well as dynamic tolls set by the system operator. Our main result shows that a broad class of decentralized congestion-dependent tolls globally stabilise the transportation network around a Wardrop equilibrium. Moreover, using dynamic marginal cost tolls, stability of the transportation network can be guaranteed around the social optimum traffic assignment.
This is particularly remarkable as the considered decentralized feedback toll policies do not require any global information about the network structure or the exogenous traffic load on the network or state and can be computed in a fully local way. We also evaluate the performance of these feedback toll policies both in the asymptotic and during the transient regime, through numerical simulations.
\end{abstract}

\section{INTRODUCTION}
In recent years. controlling the roadway congestion has become one of the main target of the transportation research community. Proposed strategies include imposing constraints on traffic flow through mechanisms such as variable speed limits, ramp metering, or traffic signal control (see \cite{HegyiSchutter1}--\cite{Varaiya} and references therein). However, such mechanisms do not consider neither the drivers' perspective nor affect the total amount of vehicles. There has been also a significant research effort to understand the drivers' answer to external communications from intelligent traveller information devices (see, e.g., \cite{Mahamassani}--\cite{Khattak}) and, in particular, studying the effect of such technologies on the drivers' route choice behaviour and on the dynamical properties of the transportation network  \cite{ComoSavla}. A traffic recommender which can announce potentially misleading travel time information and a new class of latency functions so as to influence the drivers' 
behaviour was studied in \cite{Cheng} and \cite{Amin}, respectively. Moreover, it is known that if individual drivers make their own routing decisions to minimize their own experienced delays, overall network congestion can be considerably higher than if a central planner had the ability to explicitly direct traffic. Accordingly, to charge tolls for the purpose of influencing drivers to make routing choices that result in globally optimal routing was a central research focus (see \cite{Smith}--\cite{Christ}). 

In this paper, we extend the model and results of \cite{ComoSavla} by introducing decentralized congestion-dependent tolls in order to influence the driver's route choice behaviour. Specifically, we consider a multiscale dynamical model of the transportation network whereby the traffic dynamics describing the real time evolution of the local congestion level are coupled with those of the drivers' path preferences. We assume that the latter evolve following a perturbed best response to global information about the congestion status of the whole network and to decentralized flow-dependent tolls.

Our main result shows that by using non-decreasing decentralized flow-dependent tolls and in the limit of a small update rate of the aggregate path preferences, the transportation network  globally stabilises around the Wardrop equilibrium \cite{Wardrop}. As in \cite{ComoSavla}, we assume that the drivers' path preferences evolve at a slower time scale than the physical traffic flows and adopt a singular perturbation approach \cite{Khalil} to the stability analysis of the ensuing multiscale closed-loop traffic dynamics. Note that classic results of evolutionary game theory and population dynamics \cite{Hofbauer}--\cite{SandholmLibro} cannot be applied to our framework since they suppose that the access to information take place at a single temporal and spatial scale and that the traffic dynamics are neglected by assuming that they are instantaneously equilibrated.

The introduction of tolls has long been studied as a way to influence the rational and selfish behaviour of drivers so that the associated Wardrop equilibrium can align with the system optimum network flow. A well-studied taxation mechanism that guarantees this alignment is marginal-cost pricing (see, e.g., \cite{Beckmann} and \cite{Sandholm}). Marginal-cost tolls do not require any global information about the network structure, user demands or state and can be computed in a fully local way. Using marginal-cost tolls we prove that our transportation network stabilizes around the social optimum traffic assignment. It is worth observing that our results go well beyond the traditional setting \cite{Beckmann} where only static frameworks are considered as well as \cite{Sandholm} where only path preference dynamics are consider, neglecting the physical ones that are assumed equilibrated. In fact, our analysis is carried over in a fully dynamical flow network setting. In this respect, the global optimality guarantees that are obtained in this paper through decentralized feedback toll policies should be compared with other recent results on global performance and resilience results on robust distributed control of dynamical flow networks \cite{Robust1}--\cite{Yazicioglu}.

In the last part of the paper through numerical simulations we compare the performance both asymptotic and during the transient of the system by using distributed marginal cost tolls and constant marginal cost ones. The latter, know in the literature as ``fixed" tolls (being the tolling function
on each edge a constant function of edge flow) have been well studied, and it
is known that they can be computed to enforce the social optimum equilibrium  provided that the system planner has a complete knowledge of the network topology,
user demand profile and delay functions.  
We show that not only is more convenient take into account the marginal cost tolls at convergence speed level but also they are strongly robust to variation  of network topology, user demand and traffic rate (see \cite{Brown} and \cite{BrownTAC}).

The rest of this paper is organized as follows. In Section \ref{section2}, we describe the model and observe the influence of distributed dynamics tolls on the network dynamics. In Section \ref{section3} we state the main results of the paper. 
In Section \ref{section5} we provide a numerical study of the different time and asymptotic convergences of the system.
Section \ref{section6} draws conclusions and suggests future works. Due to space limitations, we do not include any proofs of our results here and refer the reader to a forthcoming journal publication \cite{ComoMaggistro}. 
\subsection{Notation}
Let $\mathbb{R}$ and $\mathbb{R}_+ := \{x \in \mathbb{R}: x\geq 0\}$ be the set of real and nonnegative real numbers, respectively. Let $\mathcal A$ and $\mathcal B$ be finite
sets. Then $\vert \mathcal A\vert$ denotes the cardinality of $\mathcal A$, $\mathbb{R}^{\mathcal A}$ the space
of real-valued vectors whose components are indexed by
elements of $\mathcal A$, and $\mathbb{R}^{\mathcal A\times \mathcal B}$ the space of real-valued matrices whose entries are indexed by pairs in $\mathcal A\times \mathcal B$. The transpose of a matrix $Q \in \mathbb{R}^{\mathcal A\times \mathcal B}$ is denoted by $Q' \in \mathbb{R}^{\mathcal B\times \mathcal A}$, $I$ is an identity matrix and $\mathbf{1}$ an all one vector whose size depends on the context. We use the notation $\Phi:=I-\vert \mathcal A\vert^{-1}\mathbf{11'} \in \mathbb{R}^{\mathcal A\times \mathcal A}$ to denote the projection matrix of the space orthogonal to $\mathbf{1}$.
The simplex of a probability vector over $\mathcal A$ is denoted by $S(\mathcal A)=\{x \in \mathbb{R}_{+}^\mathcal{A} : \mathbf{1}'x=1\}$.
 Let $\Vert \cdot\Vert_p$ be the class of $p$-norms for $p\in [1, \infty]$, and by default, let $ \Vert \cdot \Vert:=\Vert \cdot \Vert_2$. Let now $\text{sgn}:\mathbb{R}\to \{-1, 0, 1\}$ be the sign function, defined by $\text{sgn}(x)=1$ if $x>0$, $\text{sgn}(x)=-1$ if $x<0$ and $\text{sgn}(x)=0$ if $x=0$. 
 By convention, we will assume the identity $d\vert x\vert/dx=\text{sgn}(x)$ to be valid for every $x \in \mathbb{R}$, including $x=0$. Finally, given the gradient $\nabla f$ of a function $f:D\to \mathbb{R}$ with $D\subseteq\mathbb{R}^{\mathcal A} $, we denote with  $\tilde{\nabla} f=\Phi\nabla f$  the projected gradient on $S(\mathcal A)$.
\section{MODEL DESCRIPTION}\label{section2}
\subsection{Network characteristics}

We describe the topology of the transportation network by a directed multi-graph $\mathcal{G=(V, E)}$, where $\mathcal{V}$ is a finite set of nodes
and $\mathcal{E}$ is the set of links $e$, each directed from its tail node $\theta_e$ to its head node $\kappa_e\neq \theta_e$. We shall allow for parallel links, i.e. $\theta_e=\theta_j$ and $\kappa_e=\kappa_j$ with $e\neq j$, but not for self loops, i.e., we shall assume that $\theta_e\neq \kappa_e$ for every $e \in \mathcal{E}$.
We shall denote by $B\in\{-1,0,1\}^{\mathcal V\times\mathcal E}$  the node-link incidence matrix of $\mathcal G$, whose entries are defined as $B_{ie}=1$ if $i=\theta_e$, $B_{ie}=-1$ if $i=\kappa_e$, and $B_{ie}=0$ otherwise. 
 For two nodes $o\neq d$ in $\mathcal{V}$, an $o$-$d$ path is a length-$\ell$ string of links $p=(e_1, e_2, \ldots, e_{\ell})$ such that $\theta_{e_{s+1}}=\kappa_{e_s}$ for $s=1, \ldots, \ell-1$, $\theta_{e_1}=o$, $\kappa_{e_\ell}=d$, and no node is touched twice, i.e., $i_r\neq i_s$ for all $0\leq r<s\leq \ell$. The set of $o-d$ paths in $\mathcal{G}$ of any length $\ell$ will be denoted by $\mathcal{P}$. Moreover, we shall denote the corresponding link-path incidence matrix by $A\in\{0,1\}^{\mathcal E\times\mathcal P}$ with entries 
\begin{equation*}
A_{ep}=
\begin{cases}
1&\text{ if }e\text{ is along }p\\
0&\text{otherwise}
\end{cases}
\end{equation*}
and assume that 
each link $e \in \mathcal{E}$ lies on at least one path from node $o$ to node $d$.
A path of length greater than or equal to $2$ from a node to itself is referred to as a cycle. Observe that, in contrast to \cite{ComoSavla} where the transportation network was assumed acyclic, we
allow for the presence of cycles in the network topology $\mathcal{G}$.
For every link $e \in \mathcal{E}$ and time instant $t\geq 0$ we denote the current traffic density and flow by $x_e(t)$ and $f_e(t)$ respectively, and assume the following
functional dependence
\begin{equation}\label{mufunction}
f_e=\mu_e(x_e), \quad e \in \mathcal{E},
\end{equation}
such that $\mu_e:\mathbb{R}_+\to \mathbb{R}_+$ is continuously differentiable, strictly increasing, strictly concave and 
$\mu_e(0)=0, \quad \mu_e'(0)< \infty.$ Note that in classical transportation theory the flow-density function are typically not strictly increasing, but here our
assumption is valid as long as we confine ourselves to the free-flow
region, as is done in \cite{ComoSavla}.
Then, for every link $e \in \mathcal{E}$, let
$
C_e:=\sup\{\mu_e(x_e): x_e\geq 0\}
$
be its maximum flow capacity and let 
$\mathcal{F}:=\prod_{e \in \mathcal{E}}[0, C_e)$
be the set of feasible flow vectors. 
We shall use the delay functions
$$T:\mathbb{R}_+^\mathcal{E}\to [0, +\infty]^\mathcal{E}\,,$$
\begin{equation}\label{delayfunction}
T_e(f_e):=
\begin{cases}
\displaystyle +\infty &\text{if}\ f_e\geq C_e, \\[10pt]
\displaystyle\frac{\mu_e^{-1}(f_e)}{f_e} &\text{if}\ f_e \in (0, C_e) , \\[10pt]
\displaystyle\frac{1}{\mu'_e(0)} \quad &\text{if}\ f_e=0
\end{cases}
\end{equation}
returning the delay incurred by drivers traversing link $e\in\mathcal E$, when the current flow out of it is $f_e$. 
Note that, by the properties of $\mu_e$,
$T_e(f_e)$ is continuous, strictly increasing, and such that $T_e(0) > 0$. 

\subsection{Paths choice and traffic dynamics}

We assume that the physical traffic flow consist of indistinguishable homogeneous drivers which enter in the network through the origin node, travel through it using the different paths 
and finally exit from the network through the destination node. The relative appeal of the different paths to the drivers is modelled by a time-varying probability vector over $\mathcal{P}$, which will be referred as the current \textit{aggregate path preference} and denoted by $z(t)$. Assuming a constant unit in-flow in the origin node, we consider the vector 
$$
f^{z}:=Az
$$
of the flows associated to the path preference $z(t)$ and define
$$Z:=\{z \in \mathcal{S(\mathcal{P})}:f_e^{z}<C_e\ \forall e \in \mathcal{E}\} $$ the set of feasible path preference.
The vector $z(t)$ is updated as drivers access global information about the current congestion status of the whole network (that is embodied by the flow vector $f(t)$) 
and is influenced by a vector of decentralized congestion-dependent  tolls 
\begin{equation}\label{vecTolls}
w: \mathbb{R}_+^\mathcal{E}\to [0, +\infty]^\mathcal{E}, \quad w_e(f_e)\geq 0 \quad \forall e \in \mathcal{E},
\end{equation}
that are charged to users traversing link $e$. 
In particular, we shall assume that the tolls $w_e$ are continuous and non-decreasing functions of the current flow for every link $e \in \mathcal{E}$. 

We shall assume that the cost perceived by each user crossing a link $e\in\mathcal E$ is given by the sum of the the delay $T_e(f_e)$ and the toll  $w_e(f_e)$. Moreover, as in \cite{ComoSavla}, we shall assume that path preferences are updated at some rate $\eta>0$ which is small with respect to the time scale of the network flow dynamics.
Then, from $f(t)$, the drivers evaluate the vector $A'(T(f(t))+w(f(t)))$, whose $p$th entry, $\sum_e A_{ep}(T_e(f_e(t))+w_e(f_e(t)))$, coincides with the perceived total cost that a driver expects to incur on path $p$ assuming that the congestion levels on that path won't change during the journey. Hence, according to some feasible path preference $F^{h}(f(t)) \in Z$, $z(t)$ evolves as
\begin{equation}\label{evolpi}
\dot{z}(t)=\eta(F^h(f(t))-z(t)),
\end{equation}
where $F^h:\mathcal{F}\to Z$ is a perturbed best response function, 
\begin{equation}\label{bestresponse}
F^h(f):=\operatornamewithlimits{\arg\min}_{\alpha \in Z_h}\{\alpha'A'(T(f)+w(f))+h(\alpha)\}, \quad f \in \mathcal{F},
\end{equation}
and $h:Z_h\to \mathbb{R}$ is an \textit{admissible perturbation} such that $Z_h\subseteq Z$ is a closed convex set, $h(\cdot)$ is strictly convex, twice differentiable in $int(Z_h)$, and is such that $\lim_{z\to \partial Z_h}\lVert\tilde{\nabla}h(z)\rVert=\infty$.
The definition of $F^h$ and the conditions on $h$ imply that $F^h(f) \in int(Z_h)$ and that $F^h(f)$ is differentiable on $\mathcal{F}$.\\
%
%
We now describe the \textit{local route decisions}, characterizing the fraction of drivers choosing each outgoing link when traversing a nondestination node.
Such a fraction is
the function $G_{e}(z)$ 
defined as
\begin{equation}\label{localchoice}
G_{e}(z)=
\begin{cases}
\displaystyle\frac{f_e^{z}}{\displaystyle\sum_{j \in \mathcal{E}: \theta_j=\theta_e}f_j^z} &\ \text{if} \ f_e^{z}>0, \\
\displaystyle\frac{1}{\vert\{ j \in \mathcal{E}:\theta_j=\theta_e\} \vert} &\ \text{if} \  \displaystyle\sum_{j \in \mathcal{E}: \theta_j=\theta_e}f_j^z=0,
\end{cases}
\end{equation}
for every $e \in \mathcal{E}$. Note that $\sum_k G_k(z)=1$, where $k$ are the outgoing links from the same node.\\
We refer to $
G:Z\to \mathbb{R}^\mathcal{E}$ as the local decision function that is continuously differentiable on $Z$.\\
Now, for every $e \in \mathcal{E}$
conservation of mass implies that 
\begin{equation}\label{evoldensit}
\dot{x}_e(t)=H_e(f(t), z(t)),
\end{equation}
where for all $z \in Z$ and $f \in \mathcal{F}$,
\begin{equation}\label{H}
H_e(f, z):= G_{e}(z)\bigg(\delta_{\theta_e}^{(o)}+\sum_{j: \kappa_j=\theta_e}f_j\bigg)-f_e.
\end{equation}
We now consider the evolution of the coupled dynamics
\begin{equation}\label{sistemaccoppiato}
\begin{cases}
\dot{z}(t)=\eta(F^h(f(t))-z(t)), \\
\dot{x}(t)=H(f(t), z(t))
\end{cases}
\end{equation}
where $F^h$ is defined in \eqref{bestresponse}, $\eta>0$ is the rate at which $z(t)$ is updated and $H(f,z)=\{H_e(f, z): e \in \mathcal{E}\}$.
\section{Main results}\label{section3}
In this section we give the main results of the paper.  
We shall prove that for small $\eta$ and $h$, the long-time behaviour of the system \eqref{sistemaccoppiato} is approximately at Wardrop equilibrium \cite{Wardrop} which, under proper distributed dynamic tolls, coincides with the social optimum equilibrium.
\begin{Definition}(Social optimum equilibrium). A feasible flow vector $f^* \in \mathcal{F}$ is a Social optimum equilibrium if and only if is the unique solution of the following network flow optimization problem
	\begin{equation}\label{SocialOptimum}
	f^{*}=\operatornamewithlimits{\arg\min}_{\substack{f\geq 0 \\ Bf=(\delta^{(o)}-\delta^{(d)})}}\sum_{e \in \mathcal{E}}f_eT_e(f_e).
	\end{equation}
\end{Definition}
\vspace{0.1cm}
\begin{Definition}(Wardrop equilibrium).
	For a given vector $w \in \mathbb{R}_+^\mathcal{E}$ of decentralized link tolls, a feasible flow vector $f^{(w)} \in \mathcal{F}$ is a Wardrop equilibrium if $f^{(w)}=f^z$ for some $z \in Z$ such that for all $p \in \mathcal{P}$,
	\begin{equation}\label{Wardrop}
	z_p>0 \quad \Longrightarrow \\
	\begin{array}{ll}
 \left(A'\left(T(f^z)+w(f^z)\right)\right)_p \leq\\ \left(A'\left(T(f^z)+w(f^z)\right)\right)_q \quad \forall q \in \mathcal{P}. 
	\end{array}
	\end{equation}
\end{Definition}
\vspace{0.1cm}
Existence and uniqueness of a Wardrop equilibrium are guaranteed considering the direct multi-graph $\mathcal{G}$ and under the assumption on $\mu_e$ and $w_e$.
(See Theorem 2.4 and 2.5 in \cite{Patriksson} for a complete proof). 
\begin{Theorem}\label{MainTh}
	Let $\mathcal{G}$ be the direct multi-graph, $\mu$ be as in \eqref{mufunction} and $w$ as in \eqref{vecTolls}. Then for every initial condition $(z(0), x(0)) \in Z\times [0, +\infty)^\mathcal{E}$ there exists a unique solution of \eqref{sistemaccoppiato}. Moreover, there exists a perturbed equilibrium flow $f^{(h)} \in \mathcal{F}$ such that for all $\eta>0$
	\begin{equation}
	\limsup_{t\to \infty}\lVert f(t)-f^{(h)}\rVert\leq \delta(\eta),
	\end{equation}
	where $\delta(\eta)$ is a non negative real-valued, nondecreasing function such that $\lim_{\eta\to 0}\delta(\eta)=0$. Moreover, for every sequence of admissible perturbations $\{h_k\}$ such that $\lim_{k}\lVert h_k\rVert=0$ and $\lim_{k}Z_{h_k}=\overline{Z}$ \footnote{The convergence $\lim_{k}Z_{h_k}=\overline{Z}$ holds with respect to the Hausdorff metric and $\overline{Z}$ is the closure of $Z$.}, one has
	\begin{equation}\label{limiteperturbazioni}
	\lim_{k\to \infty}f^{(h_k)}= f^{(w)}.
	\end{equation}
\end{Theorem}
\vspace{0.3cm}
Theorem \ref{MainTh} states that the system planner globally stabilises the transportation network around the Wardrop equilibrium using increasing decentralised congestion-dependent tolls.
\begin{Remark}
Note that Theorem \ref{MainTh} is not a Corollary of Theorem 2.5 in \cite{ComoSavla},
because, although the functions $T$ and $w$ both depend on the flow $f$, it is not
possible consider an auxiliary function $\overline{T} = T +w$ and directly applying the result
from \cite{ComoSavla} due to the specific structure imposed on T in \eqref{delayfunction}.
\end{Remark} 

Now, we choose as decentralized tolls the marginal cost ones, namely,
\begin{equation}\label{marginalcost}
w_e(f_e)=f_eT'_e(f_e) \quad \forall e \in \mathcal{E}.
\end{equation}
Due the properties of the delay function $T_e(f_e)$, the above tolls \eqref{marginalcost} are increasing, then the Theorem \ref{MainTh} continue to hold. Moreover the following holds 
\begin{Corollary}\label{Corollary}
Considering \eqref{marginalcost} one gets that the system \eqref{sistemaccoppiato} globally stabilises the transportation network around the social optimum traffic assignment $f^*$ without knowing arrival rates or the network structure.
\end{Corollary}
\vspace{0.1cm}
In order to prove the above we observe that considering proper costs on the links, the vector $f^{(w)}$ is the solution of a network flow optimization problem. 
Let 
$$
D_e(f_e):= \int_0^{f_e}\Big(T_e(s)+sT'_e(s)\Big)\,ds \quad e \in \mathcal{E},
$$
be the integral of the perceived cost on link $e$ using \eqref{marginalcost}. 
Then, the network flow $f^{(w)} \in \mathbb{R}_{+}^{\mathcal{E}}$ is a Wardrop equilibrium if and only if is the unique solution of the network flow optimization problem
	\begin{equation}\label{wardopcomeminimo}
	f^{(w)}=\operatornamewithlimits{\arg\min}_{\substack{f\geq 0 \\ Bf=(\delta^{(o)}-\delta^{(d)})}}\sum_{e \in \mathcal{E}}D_e(f_e),
	\end{equation}  
	where $Bf=(\delta^{(o)}-\delta^{(d)})$ is the mass conservation law. 
	Moreover, the Wardrop equilibrium coincides with the system optimum flow, 
	\begin{equation}\label{coincidenza}
	f^{(w)}=f^*.
\end{equation}
The proof of such result is very simple and use the Lagrange techniques.
\begin{Remark}\label{rem}
The tolls \eqref{marginalcost}
differ by the well now decentralized constant marginal cost tolls $w_e^{*}=f_e^{*}T'_e(f_e^{*}) \quad \forall e \in \mathcal{E}$, since the latter, in order to be used, require the knowledge of both of the social optimum flow and the inflow vector. Anyway taking into account such $w_e^{*}$, condition \eqref{coincidenza} continue to hold.
\end{Remark}

\section{Asymptotic and transient performances}\label{section5}
\begin{figure}[thpb]
	\centering
	\begin{tikzpicture}
	[scale=1.2,auto=left,every node/.style={circle,draw=black!90,scale=.5,fill=blue!30,minimum width=1cm}]
	\node (n1) at (0,0){\Large{o}};  
	\node (n2) at (2,1){\Large{a}}; 
	\node (n3) at (2,-1){\Large{b}}; 
	\node (n4) at (4,0){\Large{d}}; 
	\node [scale=0.8, auto=center,fill=none,draw=none] (n0) at (-0.8,0){};
	\node [scale=0.8, auto=center,fill=none,draw=none] (n5) at (4.8,0){};
	\foreach \from/\to in
	{n0/n1,n1/n2,n1/n3,n2/n3,n2/n4,n3/n4,n4/n5}
	\draw [-latex, right] (\from) to (\to); 
	\node [scale=2,fill=none,draw=none] (n5) at (1,0.7){$e_1$};  
	\node [scale=2,fill=none,draw=none] (n6) at (1,-0.7){$e_2$}; 
	\node [scale=2,fill=none,draw=none] (n7) at (1.85,0){$e_3$}; 
	\node [scale=2,fill=none,draw=none] (n8) at (3,0.7){$e_4$};  
	\node [scale=2,fill=none,draw=none] (n9) at (3,-0.7){$e_5$}; 
	\node [scale=1.5,fill=none,draw=none] (n10) at (-0.4,0.2){$1$};
	\node [scale=1.5,fill=none,draw=none] (n14) at (4.4, 0.2){$1$};
	\node [scale=1,fill=none,draw=none] (n11) at (4.1,0){};
	\node [scale=1,fill=none,draw=none] (n12) at (2,1.1){};
	\node [scale=1,fill=none,draw=none] (n13) at (2,-1.1){}; 	\end{tikzpicture}  
\caption{\label{graphtopology} The graph topology used for the simulations.}  
\end{figure}
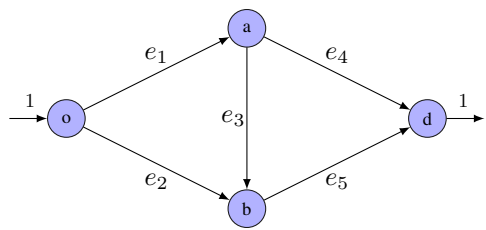

\begin{figure}[thpb]
	\centering
	\includegraphics[scale=0.55]{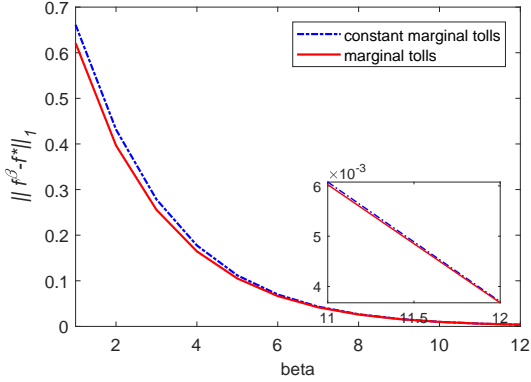}
	\caption{\label{DifferenzaDynAccopp} Plot of $\Vert f^\beta(T) - f^* \Vert_1$ for decentralised marginal and constant marginal tolls .}
\end{figure}

In this section, through numerical simulations we will compare the different performances both asymptotic and during the transient given by using the marginal cost tolls \eqref{marginalcost} and the constant marginal cost ones (see the Remark \ref{rem}).
We performed several experiments with different graph topologies for $\eta$ ranging from 0.1 to 50. In all these cases we found that the use of the decentralized marginal cost tolls is more convenient than the constant marginal ones. Indeed:\\
$-$ concerning the transient convergence, one shows that the time needed to reach the perturbed equilibrium associated to the marginal cost tolls is lower than the one to reach the equilibrium associated to the constant marginal ones;\\ 
$-$ when the admissible perturbation goes to zero, the perturbed equilibrium associated to marginal cost tolls, asymptotically converges to the social optimum flow faster than the one associated to the constant marginal cost ones.\\
We demonstrate these findings through the following example. The parameters were selected as follows:
\begin{itemize}
	\item graph topology $\mathcal{G}$ as in Fig. \ref{graphtopology};
	\item the flow-density function is 
$$
	\mu_e(x_e)=2(1-e^{-x_e}) \quad \forall e \in \mathcal{E},
$$
	and the corresponding delay function, according to \eqref{delayfunction} is given by
	\begin{equation}
	T_e(f_e)=
	\begin{cases}
	+ \infty & \text{if} \ f_e \geq  2, \\
	\displaystyle\frac{1}{f_e}\log\left(\frac{2}{2-f_e}\right) & \text{if} \ f_e\in (0, 2), \\
	1/2 & \text{if} \ f_e=0.
	\end{cases}
	\end{equation}
	\item $F^{h}$ as the logit function
\begin{equation}\label{logit}
	F_p^{h}(f)=\frac{\exp(-\beta(A'(T(f)+w(f)))_p)}{\sum_{q \in \mathcal{P}}\exp(-\beta(A'(T(f)+w(f)))_q)}, \ p \in \mathcal{P},
\end{equation}
	with $\beta>0$ the fixed noise parameter.
	\item $\eta=0.1$, $G$ as in \eqref{localchoice},
	\item initial conditions: $z_{p_1}(0)=1/2$, $z_{p_2}(0)=1/6$, $z_{p_3}(0)= 1/3$, $x_{e_1}(0)=4$, $x_{e_2}(0)=2$, $x_{e_3}(0)=3$, $x_{e_4}(0)= 1$, $x_{e_5}(0)=5$.
\end{itemize}
By the implementations follows that
for $t \in [0, 350]$ and $\beta=1$, the first time in which the system reaches the equilibrium associated to \eqref{marginalcost} is $t=2.17\cdot 10^2$, while it is $t=2.5\cdot 10^2$ the one to approach the equilibrium relative to $w_e^*$.\\
The 1-norm distance of $f^{\beta}$ 
(that is the perturbed equilibrium flow corresponding to the system \eqref{sistemaccoppiato} using \eqref{logit}) 
computed at final time $T= 350$, from the social optimum flow $f^{*}$ for $\beta$ ranging from 1 to 12 is plotted in Fig. \ref{DifferenzaDynAccopp}. This is done both considering \eqref{marginalcost} and $w_e^{*}$.
Note that the parameter $\beta$ should takes very large values in order to completely vanish the norm of the difference between $f^{\beta}$ and $f^{*}$; but, in our numerical example, we can see in Fig. \ref{DifferenzaDynAccopp} that already for $\beta=12$ the previous norm is almost null and also the asymptotic convergence of $f^\beta$ associated to \eqref{marginalcost} is slightly faster than the one of $f^\beta$ associated to $w_e^{*}$. 
\subsection{Robustness}
\begin{figure}[thpb]
	\centering
	\includegraphics[scale=0.55]{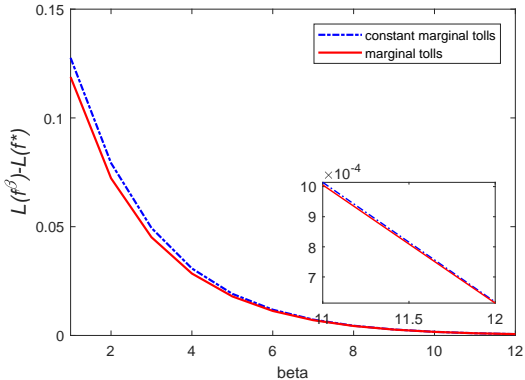}
	\caption{\label{DiffCostiDynAccoppiata} Plot of the difference $\mathcal{L}(f^{\beta}(T))-\mathcal{L}(f^{*})$ as $\beta$ increases. } 
\end{figure}

To investigate the robustness of the marginal cost tolls to variations of network's parameters, a system planner can study the effect of the variation on the total latency computed in $f^{(w)}$, where the total latency is defined as 
$$
\mathcal{L}(f)=\sum_{e \in \mathcal{E}}f_eT_e(f_e).
$$ 
By corollary \ref{Corollary} follows that the efficiency guarantees provided by the marginal cost tolls are robust to variation in network and demand structure. Indeed the following hold:
\begin{Proposition} (See \cite{Beckmann})
	For homogeneous populations, the marginal cost tolls \eqref{marginalcost} incentives optimal flows on all networks, i.e.,
	\begin{equation}\label{totallatency}
	\mathcal{L}(f^{(w)})=\mathcal{L}(f^{*}).
	\end{equation}
\end{Proposition}
\vspace{0.35cm}
Hence, the marginal cost tolls are strongly robust to variations of network topology, user demand structure and overall traffic rate.
In the following we will show (see Fig. \ref{DiffCostiDynAccoppiata}),  still using the graph topology in Fig. \ref{graphtopology} and its parameters, that
$$
\lim_{\beta\to +\infty}\mathcal{L}(f^{\beta})=\mathcal{L}(f^{*})
$$
and the asymptotic convergence using $f^\beta$ associated to \eqref{marginalcost} is lightly faster than the one in which using $f^\beta$ associated to $w_e^{*}$.
\section{CONCLUSIONS}\label{section6}
In this paper, we studied stability of Wardrop equilibria of multi-scale transportation networks with distributed dynamic tolls. We prove that if the frequency of updates of path preferences is sufficiently small and considering positive, non-decreasing decentralized flow-dependent tolls,
then the state of the network ultimately approaches a neighborhood
of the Wardrop equilibrium. Then, using a particular class of tolls, i.e., the marginal cost ones, we observe that the stability is around the social optimum equilibrium and, thanks to numerical experiments, the performances both asymptotic and during the transient of the system is better than the one obtained considering the constant marginal tolls.
In future research, inspired by the numerical results we will provide analytic estimates about the different convergence rates.
We also plan to define a more general class of tolls that does not require the knowledge of the delay functions and at the same time guarantees the convergence to the social optimum.


\begin{thebibliography}{99}
		\bibitem{HegyiSchutter1} A.~Hegyi, B.~De Schutter, and H.~Hellendoorn, ``Model predictive control for optimal coordination of ramp metering and variable speed limits," {\em Transport Res C: Emer}, vol. 13, no. 3, pp. 185-209, 2005.
	\bibitem{GomesHoro} G.~Gomes and R.~Horowitz, ``Optimal freeway ramp metering using the asymmetric cell transmission model," {\em Transport Res C: Emer}, vol. 14, no. 4, pp. 244-262, 2006.
	\bibitem{Como-convex} G.~Como, E.~Lovisari, and K.~Savla, ``Convexity and robustness of dynamic traffic assignment and freeway network control," 
{\em Transp.~Res.~B: Methodol.}, vol.~91, pp.~446--465, 2016.
\bibitem{Varaiya} P.~Varaiya, ``Max pressure control of a network of signalized intersections,'' \emph{Transport Res C: Emer.}, vol.~36, pp.~177--195, 2013.
\bibitem{Mahamassani}K.~Srinivasan and H.~Mahmassani,`` Modeling inertia and compliance mechanisms in route choice behavior under real-time information,"
{\em Transport.~Res.~Rec.}, no.~1725, pp.~45--53, 2000.
\bibitem{Khattak} A.~Khattak, A.~Polydoropoulou, and M.~Ben-Akiva, ``Modeling revealed and stated pretrip travel response to advanced traveler information
systems," {\em Transport.~Res.~Rec.}, no.~1537, pp.~46--54, 1996.
\bibitem{ComoSavla} G.~Como, K.~Savla, D.~Acemoglu, M.A.~Dahleh, and E.~Frazzoli, ``Stability analysis of transportation networks with multiscale driver decisions,"{ \em SIAM J.~Control Optim.}, vol.~51, no.~1, pp.~230--252, 2013.
 \bibitem{Cheng} Y.~Cheng and C.~Langbort, ``A model of informational nudging in transportation networks," {\em 55th IEEE Conference on Decision and Control}, pp.~7598--7604, 2016.
\bibitem{Amin} W.~Krichene, J.D.~Reilly, S.~Amin,  and A.M.~Bayen, ``Stackelberg Routing on Parallel Transportation Networks. In: Basar T., Zaccour G.~(eds) Handbook of Dynamic Game Theory. Springer, Cham, 2017.
\bibitem{Smith} M.J.~Smith, ``The marginal cost taxation of a transportation network," {\em Transp.~Res.~B: Methodol.}, vol.~13, no.~3, pp.~237--242, 1979.
\bibitem{Morrison} S.~Morrison, ``A survey of road pricing," {\em Transp.~Res.~A: Gen.}, vol.~20, no.~2, pp.~87--97, 1986.
\bibitem{Dial} R.B.~Dial, ``Network-optimized road pricing: Part I: A parable and a model," {\em Oper.~Res.}, vol.~47, pp.~54--64, 1999.
\bibitem{Cole} R.~Cole, Y.~Dodis, and T.~Roughgarden, ``How much can taxes help selfish routing?," {\em J.~Comput.~Syst.~Sci.}, vol.~72, pp.~444--467, 2006.
\bibitem{Engelson} L.~Engelson and P.~Lindberg, `` Congestion pricing of road networks with users having different time values," {\em Appl.~Optim.}, vol.~101, pp.~81-104, 2006.
\bibitem{Christ} G.~Christodoulou, K.~Mehlhorn, and E.~Pyrga, ``Improving the price of anarchy for selfish routing via coordination mechanisms ," {\em Algorithmica}, vol.~69, no.~3, pp.~619--640, 2014.
\bibitem{Wardrop}  J.G.~Wardrop, ``Some theoretical aspects of road traffic research," {\em ICE Proc.~Engrg.~Divisions}, vol.~1, no.~3, pp.~325--362, 1952.
\bibitem{Khalil} H.K.~Khalil, {\em Nonlinear Systems}, 2nd ed., Prentice-Hall, Englewood Cliffs, NJ, 1996.
\bibitem{Hofbauer}  J.~Hofbauer and K.~Sigmund, ``Evolutionary game dynamics," {\em Bull.~Amer.~Math.~Soc.}, vol.~40, pp.~479--519, 2003.
\bibitem{SandholmLibro} W.H.~Sandholm, {\em Population Games and Evolutionary Dynamics}, MIT Press, Cambridge, MA, 2011.
\bibitem{Beckmann} M.~Beckmann, C.~McGuire, and C.B.~Winsten, {\em Studies in the Economics of Transportation}, New Haven, CT: Yale University Press, 1956.
\bibitem{Sandholm} W.~Sandholm, ``Evolutionary implementation and congestion pricing," {\em Rev.~Econ.~Stud.}, vol.~69, no.~3, 667-689, 2002.


\bibitem{Robust1} G.~Como, K.~Savla, D.~Acemoglu, M.A.~Dahleh, and E.~Frazzoli. ``Robust distributed routing in dynamical networks-Part I: Locally responsive policies and weak resilience,'' {\em IEEE Trans.~Automat.~Control}, vol.~58, no.~2, pp.~317-332, 2013.
\bibitem{Robust2} G.~Como, K.~Savla, D.~Acemoglu, M.A.~Dahleh, , and E.~Frazzoli, ``Robust distributed routing in dynamical networks-Part II: Strong resilience, equilibrium selection and cascaded failures,'' {\em IEEE Trans.~Automat.~Control}, vol.~58, no.~2, 333-348, 2013.
\bibitem{Robust3} G.~Como, E.~Lovisari, and K.~Savla, `` Throughput optimality and overload behaviour of dynamical flow networks under monotone distributed routing," 
{\em IEEE Trans.~Control Netw.~Syst.}, vol.~2, no.~1, 57-67, 2015.
\bibitem{Como} G.~Como, ``On resilient control of dynamical flow networks'', \emph{Annual Reviews in Control}, vol.~43, pp.~70--80, 2017.
\bibitem{Yazicioglu} A.~Y.~Yazicioglu, M.~Roozbehani, and M.~A.~Dahleh, ``Resilient Control of Transportation Networks by Using Variable Speed Limits," {\em IEEE Trans.~Control Netw.~Syst.}, DOI 10.1109/TCNS.2017.2782364, 2017.

\bibitem{Brown} P.N.~Brown and J.R.~Marden, ``Studies on robust social influence mechanisms: Incentives for efficient network routing in uncertain settings," {\em IEEE Control Systems}, vol.~37, no.~1, pp.~98-115, 2017.
\bibitem{BrownTAC} P.N.~Brown and J.R.~Marden,``The robustness of marginal-cost taxes in affine congestion games ," {\em IEEE Trans.~Autom.~Control.}, vol.~62, no.~8, pp.~3999-4004, 2017.
\bibitem{ComoMaggistro} G.~Como and R.~Maggistro, ``On Robust Distributed Dynamic Pricing in Multiscale Transportation Networks,"  2018.


\bibitem{Patriksson} M.~Patriksson, {\em The Traffic Assignment Problem: Models and Methods}, VSP International Science, Leiden, Netherlands, 1994.
\end{thebibliography}
\end{document}